\input amstex
\documentstyle{amsppt}
\magnification=1200
\catcode`\@=11
\redefine\logo@{}
\catcode`\@=13
\pageheight{19cm}

\define \bn{\Bbb N}
\define \bz{\Bbb Z}
\define \bq{\Bbb Q}
\define \br{\Bbb R}
\define \bc{\Bbb C}

\define\Da{{\Cal D}}
\define \M{{\Cal M}}

\define\Ka{{\Cal K}}
\define \E{{\Cal E}}

\define\rk{\text{rk}~}


\define\Exc{\text{Exc}}
\define\0o{{\overline 0}}
\define\1o{{\overline 1}}
\define\a{a_1}
\redefine\b{b_1}

\define\wH{\widetilde{H}}
\define\wth{\widetilde{h}}

\TagsOnRight


\topmatter

\title
On correspondences of a K3 surface with itself, I
\endtitle

\author
Viacheslav V. Nikulin \footnote{Supported by Russian Fund of
Fundamental Research.\hfill\hfill}
\endauthor

\address
Deptm. of Pure Mathem. The University of Liverpool, Liverpool
L69 3BX, UK;
\vskip1pt
Steklov Mathematical Institute,
ul. Gubkina 8, Moscow 117966, GSP-1, Russia
\endaddress
\email
vnikulin\@liv.ac.uk\ \
vvnikulin\@list.ru
\endemail

\abstract
Let $X$ be a K3 surface with a polarization $H$ of degree
$H^2=2rs$, $r,\,s\ge 1$.
Assume $H\cdot N(X)=\bz$ for the Picard lattice
$N(X)$. The moduli space of sheaves over $X$ with the isotropic Mukai
vector $(r,H,s)$ is again a K3 surface $Y$.

We prove that $Y\cong X$, if there exists
$h_1\in N(X)$ with $(h_1)^2=f(r,s)$, $H\cdot h_1\equiv 0\mod g(r,s)$,
and $h_1$ satisfies some condition of primitivity. These conditions
are necessary, if $X$ is general with $\rk N(X)=2$. Existence of such kind
a criterion is surprising, and it also gives some geometric
interpretation of elements in $N(X)$ with negative square.

We describe all irreducible 18-dimensional components of moduli of
the $(X,H)$ with $Y\cong X$. We prove that their number is always infinite.

This generalizes results of \cite{4, 5} for $r=s$.
\endabstract

\rightheadtext
{On correspondences of K3 with itself}
\leftheadtext{V.V. Nikulin}
\endtopmatter

\document

\head
0. Introduction
\endhead

Let $X$ be a K3 surface with a polarization $H$ of degree $H^2=2rs$,
$r,\,s\ge 1$, and a primitive Mukai vector $(r,H,s)$. Let $Y$ be the
moduli of sheaves (coherent) over $X$ with the Mukai
vector $(r,H,s)$.
The Mukai vector $v=(r,H,s)$ is isotropic, then $Y$ is a K3 surface
with a natural $nef$ element $h$ with $h^2=2ab$
where $c=\text{g.c.d}(r,s)$ and $a=r/c$, $b=s/c$.
The surface $Y$ is isogenous to $X$ in the sense of Mukai.
The second Chern class of the corresponding quasi-universal sheave gives
a 2-dimensional algebraic cycle $Z\subset X\times Y$, and an
algebraic correspondence between $X$ and $Y$. See \cite{6, 7}
and also \cite{1} about these results.

Let $N(X)$ be the Picard lattice of $X$. We consider an
invariant $\gamma (H)\in \bn$ of $H\in N(X)$ defined by
$H\cdot N(X)=\gamma (H) \bz$ (clearly, $\gamma (H)|2rs$). In this
paper, we assume that
$$
\gamma (H)=1,\ \text{that is}\ \ H\cdot N(X)=\bz .
\tag{0.1}
$$
Then $H$ is primitive and $\rho (X)=\rk N(X)\ge 2$.
By Mukai \cite{6, 7}, the transcendental
periods $(T(X), H^{2,0}(X))$ and $(T(Y),H^{2,0}(Y))$ are
isomorphic in this case. We can expect that surfaces $X$
and $Y$ also can be isomorphic, and we then get a cycle
$Z\subset X\times X$, and a correspondence of $X$ with itself.
Thus, an interesting for us question is

\proclaim{Question 1} When is $Y$ isomorphic to $X$?
\endproclaim

We want to answer this question in terms of Picard lattices
$N(X)$ and $N(Y)$ of $X$ and $Y$. Then our question can be
reformulated as follows:

\proclaim{Question 2} Assume that $N$ is a hyperbolic lattice,
$\widetilde{H}\in N$ an element with square $2rs$ and
$\widetilde{H}\cdot N=\bz$. What are
conditions on $N$ and $\widetilde{H}$ such that for any K3 surface
$X$ with Picard lattice $N(X)$ and s polarization $H\in N(X)$
the corresponding K3 surface $Y$ is isomorphic to $X$,
if the pairs $(N(X), H)$ and $(N, \widetilde{H})$ are
isomorphic as abstract lattices with fixed elements?

In other words, what are conditions on $(N(X), H)$ as an abstract lattice
with a vector $H$ which are sufficient for
$Y$ to be isomorphic to $X$ and are necessary,
if $X$ is a general K3 surface with its Picard lattice $N(X)$?
\endproclaim

We answered this question in \cite{4} for $r=s=2$ (then $H^2=8$,
and the condition \thetag{0.1} always satisfies,
if $H$ is primitive); in \cite{5} for $r=s$ under the condition
\thetag{0.1} (then \thetag{0.1} always satisfies, if $H$ is
primitive and $r=s$ is even). The main surprising result of
\cite{4, 5}  was that
$Y\cong X$, if the Picard lattice $N(X)$ has an element $h_1$ with some
prescribed square $(h_1)^2$ and some minor additional conditions,
and these conditions are necessary for a general K3 surface
$X$ with $\rho (X)=\rk N(X)=2$. Here we prove similar results in general,
for arbitrary $r$ and $s$, under the condition \thetag{0.1}.

\proclaim{Theorem 1}
Let $X$ be a K3 surface with a polarization $H$ such that
$H^2=2rs$, $r,\,s\ge 1$, the Mukai vector $(r,H,s)$ is primitive
and $Y$ the moduli of sheaves on $X$ with the Mukai vector
$(r,H,s)$.

Then $Y\cong X$, if at least for one of signs $\pm$ there exists
$h_1\in N(X)$ such that the elements $H$, $h_1$ are contained in a
2-dimensional sublattice $N\subset N(X)$ with $H\cdot N=\bz$, and $h_1$
belongs to one of $a$-series or $b$-series described below:

$a$-series:
$$
h_1^2=\pm 2bc,\  H\cdot h_1\equiv 0\mod bc,\
H\cdot h_1\not\equiv 0\mod bcl_1,\
h_1/l_2\not\in  N(X)
$$
where $l_1$ and $l_2$ are any primes such that $l_1^2|a$ and $l_2^2|b$;

$b$-series:
$$
h_1^2=\pm 2ac,\  H\cdot h_1\equiv 0\mod ac,\
H\cdot h_1\not\equiv 0\mod acl_1,\
h_1/l_2\not\in  N(X).
$$
where $l_1$ and $l_2$ are any primes such that $l_1^2|b$ and $l_2^2|a$.

The conditions above are necessary for $\gamma(H)=1$ and
$Y\cong X$, if $\rho(X)\le 2$ and $X$ is a general K3 surface with
its Picard lattice, i. e. the automorphism group of the transcendental
periods $(T(X), H^{2,0}(X))$ of $X$ is $\pm 1$.
\endproclaim

Like in \cite{4, 5} for $r=s$,
we also describe all (irreducible) divisorial conditions on
19-dimensional moduli of polarized K3 surfaces $(X,H)$
which imply $\gamma (H)=1$ and
$Y\cong X$. It is sufficient to assume that $\rho (X)=2$.
According to \cite{10, 11},
it is equivalent to description up to isomorphisms of all possible pairs
$H\in N(X)$ where $\rk N(X)=2$.
There are two defining invariants of such pairs:
$d=-\det N(X)\in \bn$ and
$\overline{\mu}=\{\mu,-\mu\}\subset (\bz/2rs)^\ast$
where $d\equiv \mu^2\mod 4rs$.
See the definition of $\mu$ in \thetag{3.1.2} and \thetag{3.1.3}).

We introduce sets:

The $a$-series:
$$
\align
\Da(r,s;a)^{\overline{\mu}}_\pm
=\{
& d\in \bn \ | \
d\equiv \mu^2 \mod 4abc^2, \
\exists \
(p,q) \in \bz \times \bz:\  p^2-dq^2=\pm 4ac,\  \\
&p-\mu q \equiv 0 \mod 2ac,\
\text{g.c.d}(a,p,q)=1,\\
&p-\mu q\not\equiv 0 \mod 2ac l\
\text{for any prime}\ l\ \text{such that}\ l^2|b\ \}.
\endalign
$$

The $b$-series:
$$
\align
\Da(r,s;b)^{\overline{\mu}}_\pm
=\{
& d\in \bn \ | \
d\equiv \mu^2 \mod 4abc^2, \
\exists \
(p,q) \in \bz \times \bz:\  p^2-dq^2=\pm 4bc,\  \\
&p-\mu q \equiv 0 \mod 2bc,\
\text{g.c.d}(b,p,q)=1,\\
&p-\mu q\not\equiv 0 \mod 2bc l\
\text{for any prime}\ l\ \text{such that}\ l^2|a\ \}.
\endalign
$$

We have

\proclaim{Theorem 2}
With the above notations all possible irreducible divisorial conditions on
moduli of polarized K3 surfaces $(X,H)$ with a polarization
$H$ with $H^2=2rs$, which imply $\gamma (H)=1$ and
$Y\cong X$ are labelled by the set
$$
{\Cal Div}(r,s)={\Cal Div}(r,s;a)\bigcup {\Cal Div}(r,s;b),
$$
where
$$
\split
{\Cal Div}(r,s;a)=\{(d,\,\overline{\mu})\ |\
\overline{\mu}=\{\mu, -\mu\}\subset
(\bz/2rs)^\ast,\ d\in \Da(r,s;a)^{\overline{\mu}}\},\\
{\Cal Div}(r,s;b)=\{(d,\,\overline{\mu})\ |\
\overline{\mu}=\{\mu, -\mu\}\subset
(\bz/2rs)^\ast,\ d\in \Da(r,s;b)^{\overline{\mu}}\}
\endsplit
$$
and
$$
\Da(r,s;a)^{\overline{\mu}}=\Da(r,s;a)^{\overline{\mu}}_+
\bigcup \Da(r,s;a)^{\overline{\mu}}_-,\ \
\Da(r,s;b)^{\overline{\mu}}=\Da(r,s;b)^{\overline{\mu}}_+
\bigcup \Da(r,s;b)^{\overline{\mu}}_-.
$$

The sets
$$
\Da(r,s;a)^{ \overline{\mu}}_{\pm}\ \text{and}\
\Da(r,s;b)^{ \overline{\mu}}_{\pm}
$$
are infinite if they are not empty.

The set ${\Cal Div}(r,s)$ is always non-empty and then infinite: i. e.
at least one of sets $\Da(r,s;a)^{ \overline{\mu}}_{\pm}$ or
$\Da(r,s;b)^{ \overline{\mu}}_{\pm}$ is not empty for one of
$\overline{\mu}$ and $\pm$.

See more exact and strong statements in Lemmas 4.3, 4.4 and
Theorems 4.1, 4.5 and 4.6.
\endproclaim

These two theorems give the main results of the paper. We remark that
Theorem 1 gives many (if not all) known cases when $Y\cong X$ (e.g.
see \cite{2}, \cite{17---19} and \cite{20}). Theorem 1 is
also interesting because it gives a very clear geometric interpretation
of elements of $N(X)$ with negative square (when the sign is $-$).
For K3 surfaces it is well-known only for elements $\delta\in N(X)$ with
$\delta^2=-2$: then $\delta$ or $-\delta$ is effective.

In Section 5 we formulate similar results for the surface $Y$.
They are very similar to the results for $X$.
We only formulate these results since proofs are very similar
to the corresponding proofs for $X$, but
sometimes require additional non-trivial considerations.
Thus, it is surprising that the results for $Y$ are so similar to the
results for $X$. Using results for $X$ and $Y$ together, one can
construct non-trivial compositions of correspondences of the
K3 surfaces $X\cong Y$ with itself.

In ideas   we always follow \cite{4} where a particular case
$r=s=2$ had been treated in great details. But technically the
results of this paper are much more complicated. One needs to analyze
many details and calculations in general and find appropriate
formulations. As in \cite{4, 5}, the fundamental tools to get
the results above is the Global Torelli Theorem for K3 surfaces
proved by I.I.Piatetsky-Shapiro and I.R.Shafarevich in \cite{14},
and results of Mukai \cite{6, 7}.
By results of \cite{6, 7},
we can calculate periods of $Y$ using periods of $X$;
by the Global Torelli Theorem for K3 surfaces \cite{14},
we can find out if $Y$ is isomorphic to $X$.

Recently the general case of a primitive isotropic Mukai vector
$(r,H,s)$ was considered in \cite{13} where many results of this paper
were generalized to the general case when $\gamma(H)$ is possibly
not equal to one and $H$ is possibly not primitive in $N(X)$.
Results of \cite{13} are much more complicated and much less efficient.
It is non-trivial to deduce results of this paper
directly from results of \cite{13}.

Thus, in this paper we treat a very nice and important particular
case when the results are very simple, efficient and complete.
Using methods of this paper and \cite{13}, one can extend
these results to more general types of primitive isotropic Mukai
vectors $(r,H,s)$. We hope to consider that in further publications.

I am grateful to Carlo Madonna for useful discussions.

\head
1. Some notations and results about lattices and K3 surfaces
\endhead

\subhead
1.1. Some notations about lattices
\endsubhead
We use notations and terminology from \cite{11} about lattices,
their discriminant groups and forms. A {\it lattice} $L$ is a
non-degenerate integral symmetric bilinear form. That is $L$ is a free
$\bz$-module equipped with a symmetric pairing $x\cdot y\in \bz$ for
$x,\,y\in L$, and this pairing should be non-degenerate. We denote
$x^2=x\cdot x$. The {\it signature} of $L$ is the signature of the
corresponding real form $L\otimes \br$.
The lattice $L$ is called {\it even}
if $x^2$ is even for any $x\in L$. Otherwise, $L$ is called {\it odd}.
The {\it determinant} of $L$ is defined to be $\det L=\det(e_i\cdot e_j)$
where $\{e_i\}$ is some basis of $L$. The lattice $L$ is {\it unimodular}
if $\det L=\pm 1$.
The {\it dual lattice} of $L$ is
$L^\ast=Hom(L,\,\bz)\subset L\otimes \bq$. The
{\it discriminant group} of $L$ is $A_L=L^\ast/L$. It has the order
$|\det L|$. The group $A_L$ is equipped with the
{\it discriminant bilinear form} $b_L:A_L\times A_L\to \bq/\bz$
and the {\it discriminant quadratic form} $q_L:A_L\to \bq/2\bz$
if $L$ is even. To get this forms, one should extend the form of $L$ to
the form on the dual lattice $L^\ast$ with values in $\bq$.

For $x\in L$, we shall consider
the invariant $\gamma (x)\ge 0$ where
$$
x\cdot L=\gamma (x)\bz .
\tag{1.1.1}
$$
Clearly, $\gamma (x)|x^2$ if $x\not=0$.

We denote by $L(k)$ the lattice obtained from a lattice $L$ by
multiplication of the form of $L$ by $k\in \bq$.
The orthogonal sum of lattices $L_1$ and $L_2$ is denoted by
$L_1\oplus L_2$.
For a symmetric integral matrix
$A$, we denote by $\langle A \rangle$ a lattice which is given by
the matrix $A$ in some bases. We denote
$$
U=\left(
\matrix
0&1\\
1&0
\endmatrix
\right).
\tag{1.1.2}
$$
Any even unimodular lattice of the signature $(1,1)$ is isomorphic to
$U$.

An embedding $L_1\subset L_2$ of lattices is called {\it primitive}
if $L_2/L_1$ has no torsion.
We denote by $O(L)$, $O(b_L)$ and $O(q_L)$ the automorphism groups of
the corresponding forms.
Any $\delta\in L$ with $\delta^2=-2$ defines
a reflection $s_\delta\in O(L)$ which is given by the formula
$$
x\to x+(x\cdot \delta)\delta,
$$
$x\in L$. All such reflections generate
the {\it 2-reflection group} $W^{(-2)}(L)\subset O(L)$.

\subhead
1.2. Some notations about K3 surfaces
\endsubhead
Here we remind some basic notions and results about K3 surfaces,
e. g. see \cite{14---16}.
A K3 surface $X$ is a non-singular projective algebraic surface over
$\bc$ such that its canonical class $K_X$ is zero and the irregularity
$q_X=0$. We denote by $N(X)$ the {\it Picard lattice} of $X$ which is
a hyperbolic lattice with the intersection pairing
$x\cdot y$ for $x,\,y\in N(S)$. Since the canonical class $K_X=0$,
the space $H^{2,0}(X)$ of 2-dimensional holomorphic differential
forms on $X$ has dimension one over $\bc$, and
$$
N(X)=\{x\in H^2(X,\bz)\ |\ x\cdot H^{2,0}(X)=0\}
\tag{1.2.1}
$$
where $H^2(X,\bz)$ with the intersection pairing is a
22-dimensional even unimodular lattice of signature
$(3,19)$. The orthogonal lattice $T(X)$ to $N(X)$ in $H^2(X,\bz)$ is called
the {\it transcendental lattice of $X$.} We have
$H^{2,0}(X)\subset T(X)\otimes \bc$. The pair $(T(X), H^{2,0}(X))$ is
called the {\it transcendental periods of $X$}.
The {\it Picard number} of $X$ is
$\rho(X)=\rk N(X)$. A non-zero element $x\in N(X)\otimes \br$ is
called {\it nef} if $x\not=0$ and $x\cdot C\ge 0$ for any effective
curve $C\subset X$. It is known that an element $x\in N(X)$ is ample
if $x^2>0$, $x$ is $nef$, and the orthogonal complement
$x^\perp$ to $x$ in  $N(X)$ has no elements with square $-2$.
For any element $x\in N(X)$ with $x^2\ge 0$, there exists an element
$w\in W^{(-2)}(N(X))$ such that
$\pm w(x)$ is nef; it then is ample, if $x^2>0$ and $x^\perp$ had no
elements with square $-2$ in $N(X)$.

We denote by $V^+(X)$ the light cone of $X$, which is the half-cone of
$$
V(X)=\{x\in N(X)\otimes \br\ |\ x^2>0\ \}
\tag{1.2.2}
$$
containing a polarization of $X$. In particular, all $nef$ elements
$x$ of $X$ belong to $\overline{V^+(X)}$:
one has $x\cdot V^+(X)>0$ for them.

The reflection group $W^{(-2)}(N(X))$ acts in $V^+(X)$ discretely,
and its  fundamental
chamber is the closure $\overline{\Ka(X)}$ of the K\"ahler cone
$\Ka(X)$ of $X$. It is the same as the set of all $nef$ elements of $X$.
Its faces are orthogonal to the set $\Exc(X)$ of all exceptional curves $r$
on $X$ which are non-singular rational curves $r$ on $X$ with $r^2=-2$.
Thus, we have
$$
\overline{\Ka (X)}=\{0\not=x\in \overline{V^+(X)}\ |
\ x\cdot \Exc(X)\ge 0\,\}.
\tag{1.2.3}
$$
\par\bigskip

\head
2. General results on the Mukai correspondence between K3 surfaces
with primitive polarizations of degrees $2rs$ and $2ab$
which gives isomorphic K3's
\endhead

\subhead
2.1. The correspondence
\endsubhead
Let $X$ be a smooth complex projective K3 surface with a
polarization $H$ of degree $2rs$, $r,\,s\ge 1$.
Let $Y$ be the moduli space of coherent sheaves $\E$ on $X$ with the
primitive isotropic Mukai vector $v=(r,H,s)$.
Then $\rk \E=r$, $\chi (\E)=r+s$ and $c_1(\E)=H$. Let
$$
H^\ast(X,\bz)=H^0(X,\bz)\oplus H^2(X,\bz)\oplus H^4(X,\bz)
\tag{2.1.1}
$$
be the full cohomology lattice of $X$ equipped with the Mukai product
$$
(u,v)=-(u_0\cdot v_2+u_2\cdot v_0)+u_1\cdot v_1
\tag{2.1.2}
$$
for $u_0,v_0\in H^0(X,\bz)$, $u_1,v_1\in H^2(X,\bz)$,
$u_2,v_2\in H^4(X,\bz)$. We naturally identify
$H^0(X,\bz)$ and $H^4(X,\bz)$ with $\bz$. Then the Mukai product is
$$
(u,v)=-(u_0v_2+u_2v_0)+u_1\cdot v_1.
\tag{2.1.3}
$$
Since $H^2=2rs$, the element
$$
v=(r,H,s)=(r,H,\chi-r)\in H^\ast(X,\bz)
\tag{2.1.4}
$$
is isotropic, i.e. $v^2=0$. In this case (for a primitive $v$),
Mukai showed \cite{6, 7} that $Y$ is a K3 surface, and
one has the natural identification
$$
H^2(Y,\,\bz) \cong (v^\perp/\bz v)
\tag{2.1.5}
$$
which also gives the isomorphism of the Hodge structures of $X$ and $Y$.
The $Y$ has the canonical polarization $h=(-a,0,b)\mod \bz v$ where
$$
c=\text{g.c.d}(r,s),\ a=r/c\,,\ b=s/c\,.
\tag{2.1.6}
$$
Then $\text{g.c.d.}(a,b)=1$.

In particular, \thetag{2.1.5} gives an embedding
$$
T(X)\subset T(Y)
\tag{2.1.7}
$$
of the lattices of transcendental cycles of the index
$$
[T(Y):T(X)]=n(v)=\min {|v\cdot x|}
\tag{2.1.8}
$$
where $x\in H^0(X,\bz)\oplus N(X)\oplus H^4(X,\bz)$ and $v\cdot x\not=0$.

In this paper, we are interested in the case when
$Y\cong X$. By \thetag{2.1.8}, it may happen, if $n(v)=1$ only.

We can introduce the invariant $\gamma=\gamma (H)\in \bn$ which
is defined by
$$
H\cdot N(X)=\gamma \bz.
\tag{2.1.9}
$$
In this paper we consider the important case when $\gamma=\gamma (H)=1$.
Thus, $N(X)\cdot H=\bz$. Equivalently, there exists $x\in N(X)$ such that
$x\cdot H=1$.
From $\gamma=\gamma (H)=1$, it follows that
$H\in N(X)$ is primitive, the Mukai
vector $(r,H,s)$ is also primitive, and the invariant
$n(v)=1$. Thus, we may have a hope to have cases when $Y\cong X$.
Since the Picard lattice $N(X)$ is even, it also follows that
the Picard number $\rho(X)=\rk N(X)\ge 2$ in this case. Thus:
$$
\split
&\gamma =\gamma (H)=1 \ \Longrightarrow\\
&\rho(X)\ge 2,\
\text{$H$ is primitive},\ (T(Y),H^{2,0}(Y))\cong (T(X), H^{2,0}(X)).
\endsplit
\tag{2.1.10}
$$

The same case $\gamma =1$ was considered in \cite{4} for $r=s=c=2$,
and in \cite{5} for $r=s=c$. Results of \cite{4, 5}
concerning an arbitrary Picard lattice $N(X)$ can be generalized to
the general case of $\gamma =1$.

\subhead
2.2. General results for an arbitrary Picard lattice $N(X)$
\endsubhead
First we remind {\it the characteristic map} defined in \cite{4}.
It is defined for any non-degenerate even lattice $S$ and any
its element $P\in S$ with $P^2=2m\not=0$ and $\gamma (P)=1$, i. e.
$P\cdot S=1$. Then $P$ is primitive in $S$.

We denote by
$$
K(P)=P^\perp_{S}
\tag{2.2.1}
$$
the orthogonal complement to $P$ in $S$.
Put $P^\ast=P/2m$. Then any element $x\in S$ can be written as
$$
x=nP^\ast+k^\ast
\tag{2.2.2}
$$
where $n\in \bz$ and $k^\ast\in K(P)^\ast$, because
$$
\bz P\oplus K(P)\subset S\subset
S^\ast\subset \bz P^\ast\oplus K(P)^\ast.
$$
Since $\gamma(P)=1$,
the map $nP^\ast+[P] \to k^\ast +K(P)$ gives an isomorphism of
the groups
$\bz/2m \cong [P^\ast]/[P]\cong [u^\ast(P)+K(P)]/K(P)$ where
$u^\ast(P)+K(P)$ has order $2m$ in $A_{K(P)}=K(P)^\ast/K(P)$.
It follows,
$$
S=[\bz P, K(P), P^\ast+u^\ast(P)].
\tag{2.2.3}
$$
The element $u^\ast(P)$ is defined canonically mod $K(P)$. Since
$P^\ast+u^\ast(P)$ belongs to the even lattice $S$, it follows
$$
(P^\ast+u^\ast(P))^2=\frac{1}{2m}
		   +{u^\ast(P)}^2 \equiv 0 \mod 2.
\tag{2.2.4}
$$
Let $\overline{P^\ast}=P^\ast \mod [P]\in [P^\ast]/[P]\cong \bz/2m$
and $\overline{k^\ast}=k^\ast\mod K(P)\in A_{K(P)}$.
Then
$$
S/[P,K(P)]=(\bz/2m)(\overline{P^\ast} +\overline{u^\ast(P)})\subset
(\bz/2m)\overline{P^\ast}+K(P)^\ast/K(P).
\tag{2.2.5}
$$
Also $S^\ast\subset \bz P^\ast+K(P)^\ast$ since
$P+K(P)\subset S$,
and for $n\in \bz$, $k^\ast \in K(P)^\ast$ we have
$x=nP^\ast+k^\ast \in S^*$, if and only if
$$
(nP^\ast+k^\ast)\cdot (P^\ast+u^\ast(P))={n\over 2m}+
k^\ast\cdot u^\ast(P)\in \bz.
$$
It follows,
$$
\aligned
S^\ast &=
\{nP^\ast+k^\ast\ |\ n\in \bz,\ k^\ast \in K(P)^\ast,\
n\equiv -2m \ (k^\ast\cdot u^\ast(P)) \mod 2m \}\subset \\
& \subset \bz P^\ast+K(P)^\ast,
\endaligned
\tag{2.2.6}
$$
and
$$
S^\ast/[P,K(P)]=
\{-2m(\overline{k^\ast}\cdot \overline{u^\ast(P)})\, \overline{P^\ast}+
\overline{k^\ast}\}\ |\ \overline{k^\ast} \in A_{K(P)}\}
\subset (\bz/2m) \overline{P^\ast}+A_{K(P)}.
\tag{2.2.7}
$$
We introduce {\it the characteristic map of $P\in S$:}
$$
\kappa(P):K(P)^\ast \to A_{K(P)}/(\bz/2m)(u^\ast(P)+K(P))\to A_{S}
\tag{2.2.8}
$$
where for $k^\ast \in K(P)^\ast$ we have
$$
\kappa(P)(k^\ast)=-2m(k^\ast\cdot u^\ast(P)) P^\ast+k^\ast + S\in
A_{S}.
\tag{2.2.9}
$$
It is epimorphic and its kernel is $(\bz/2m)(u^\ast(P)+K(P))$. It gives
the canonical isomorphism
$$
\overline{\kappa(P)}:A_{K(P)}/(\bz/2m)(u^\ast(P)+K(P))
\cong A_{S}.
\tag{2.2.10}
$$
For the corresponding discriminant forms we have
$$
\kappa(k^\ast)^2 \mod 2=(k^\ast)^2+2m(k^\ast\cdot u^\ast(P))^2\mod 2.
\tag{2.2.11}
$$
Below we shall consider the characteristic map for elements
$P\in N(X)$ with $\gamma (P)=1$. It will be denoted as $\kappa (P)$.

We define $m(a,b)\mod 2ab$ by
$$
m(a,b)\equiv -1\mod 2a,\ \ \text{and}\ \  m(a,b)\equiv 1\mod 2b.
\tag{2.2.12}
$$
The $m(a,b)\mod 2ab$ was introduced by Mukai \cite{8}.

The same calculations as in \cite{4} and also \cite{5}
(the key calculations were done in \cite{4}) prove the fundamental
for further results Theorem 2.2.1 below. Its complete proof in the
general case (for any $\gamma$ and possibly non-primitive $H$)
is contained in \cite{13}. The Global Torelli Theorem for K3 surfaces
\cite{14} and Mukai's result \thetag{2.1.5} are fundamental for
the proof.

\proclaim{Theorem 2.2.1} Assume that $X$ is a K3 surface
with a polarization $H$ with $H^2=2rs$, $r,\,s\ge 1$, and $\gamma (H)=1$
for  $H\in N(X)$. Let $Y$ be the moduli space of coherent sheaves on
$X$ with the Mukai vector $v=(r,H,s)$ (it is clearly primitive).

Then the surface $Y$ is isomorphic to $X$,
if the following conditions (i) and (ii) hold:

(i) there exists $\widetilde{h}\in N(X)$ of the degree
$\widetilde{h}^2=2ab$, $\gamma(\widetilde{h})=1$, and such that
there exists an embedding
$$
f:K(H)\rightarrow K(\wth)
$$
of negative definite lattices such that
$$
\split
&f^\ast (K(\wth))=[K(H),2abcu^\ast (H)],\\
&f^\ast(u^\ast (\wth))+f^\ast(K(\wth))=m(a,b)c\,u^\ast(H)+f^\ast(K(\wth)).
\endsplit
\tag{2.2.13}
$$

(ii) For the $\wth$  and $f$ in (i), there exists a choice of
$\pm$ such that
$$
\kappa(\wth)(k^\ast)=\pm \kappa(H)(f^\ast(k^\ast))
\tag{2.2.14}
$$
for any $k^\ast \in K(\widetilde{h})^\ast$.

The conditions (i) and (ii) are necessary for a K3 surface $X$ with
$\rho (X)\le 19$ which is general for its Picard lattice $N(X)$
in the following sense: the automorphism group of the transcendental periods
$(T(X),H^{2,0}(X))$ is $\pm 1$. If $\rho (X)=20$, then always $Y\cong X$.
\endproclaim

For instance, let $c=1$. Then one can take $\wth=H$ and $f=id$
in Theorem 2.2.1, if $m(a,b)\equiv 1\mod 2ab$. This is valid, only if either
$a=1$ or $b=1$. Thus, $Y\cong X$, if $\gamma (H)=1$, $c=1$ and one of $a$,
$b$ is equal to $1$. According to Mukai \cite{8}, always $Y\cong X$, if
$c=1$ and one of $a$, $b$ is equal to $1$. One does not need the
additional condition $\gamma (H)=1$.

\head
3. The case of Picard number 2
\endhead

\subhead
3.1. Main results for $\rho(X)=2$
\endsubhead
Here we apply results of Sect. 2 to $X$ and $Y$ with
Picard number 2. A general K3 surface $X$ with $\gamma(H)=1$ has
$\rho (X)=2$.

First, we consider an arbitrary
K3 surface $X$ with Picard number 2 and a primitive polarization $H$ of
degree $H^2=2rs$ where $r,\,s\ge 1$. Like above, $c=\text{g.c.d.}(r,s)$ and
$r=ca$, $s=cb$.
Additionally, we assume that $\gamma (H)=1$ for $H\in N(X)$.
Let
$$
K(H)=H^\perp_{N(X)}=\bz \delta
\tag{3.1.1}
$$
and $\delta^2=-t$ where $t>0$ is even. The $\delta\in N(X)$ is defined
uniquely up to $\pm \delta$. Since $\gamma (H)=1$, we have
$$
N(X)=[\bz H,\, \bz \delta,\, \mu H^\ast+{\delta \over 2rs}]
\tag{3.1.2}
$$
where $H^\ast=H/2rs$ and $g.c.d(\mu,2rs)=1$.
The element
$$
\pm \mu \mod 2rs \in   (\bz/2rs)^\ast
\tag{3.1.3}
$$
is {\it the invariant of the pair $(N(X),H)$}
up to isomorphisms of lattices
with a primitive vector $H$ of $H^2=2rs$ and $\gamma (H)=1$.
If $\delta $ changes to
$-\delta$, then $\mu\mod 2rs$ changes to $-\mu\mod 2rs$.
We have
$$
(\mu H^\ast+{\delta\over 2rs})^2=
\frac{1}{2rs}(\mu^2-\frac{t}{2rs}) \equiv 0\ \mod 2.
\tag{3.1.4}
$$
Then $t=2rs d$, for some $d \in \bn$ and
$\mu^2 \equiv\ d \mod 4rs$. Thus, $d\mod 4rs \in (\bz/4rs)^{\ast\,2}$.
We have $-d=\det (N(X))$.
Any element $z\in N(X)$ can be written as
$z=(xH+y\delta)/2rs$ where $x \equiv \mu y \mod 2rs$.
In these considerations, one can replace $H$ by
any primitive element of $N(X)$ with square $2rs$ and
$\gamma (H)=1$. Thus, we have:

\proclaim{Proposition 3.1.1} Let $X$ be a K3 surface with Picard
number $\rho = 2$ equipped with a primitive
polarization $H$ of degree $H^2=2rs$, $r,s\ge 1$, and
$\gamma (H)=1$.

The pair $(N(X),H)$ has the invariants $d\in \bn$ and
$\pm \mu\mod 2rs\in (\bz/2rs)^\ast$
such that $\mu^2\equiv d \mod 4rs$. We have: $\det N(X)=-d$, and
$K(H)=H^\perp_{N(X)}=\bz \delta$ where $\delta^2=-2rsd$. Moreover,
$$
N(X)=[H,\delta,\frac{(\mu H+\delta)}{2rs}],
\tag{3.1.5}
$$
$$
N(X)=\{z=\frac{(xH+y\delta)}{2rs} \ |\ x,y\in \bz\ \text{and}\
x \equiv \mu y \mod 2rs \}.
\tag{3.1.6}
$$
We have $z^2=(x^2-dy^2)/2rs$.

For any primitive element $H^\prime \in N(X)$ with
$(H^\prime)^2 =H^2=2rs$ and the same invariant $\pm \mu$,
there exists an automorphism
$\phi\in O(N(X))$ such that $\phi(H)=H^\prime$.
\endproclaim

From Proposition 3.1.1 and $\text{g.c.d}(d,2rs)=\text{g.c.d}(d,2ab)=1$
where $\det N(X)=-d$, we get

\proclaim{Proposition 3.1.2}
Under conditions and notations of
Propositions 3.1.1, all elements
$h^\prime=(xH+y\delta)/(2rs) \in N(X)$ with
$(h^\prime)^2=2ab$ are in one to one
correspondence with integral solutions $(x,y)$ of the equation
$$
x^2-dy^2=4a^2b^2c^2
\tag{3.1.7}
$$
such that $x \equiv \mu y \mod 2abc^2$.

For $k\in \bn$ and $k>1$ the element $h^\prime$ is divisible by $k$
in $N(X)$ and then it is not primitive, if and only if $k\,|\, (x,y)$ and
$x\equiv \mu y\mod 2abc^2 k$. Then $k^2|ab$. We call such a pair
$(x,y)$ as $\mu$-divisible by $k$. If $(x,y)$ is not $\mu$-divisible by
$k>1$, we call such a pair $(x,y)$ as $\mu$-primitive.
It is enough to consider only prime $k$. Thus, $h^\prime$ is primitive,
if and only if the corresponding pair $(x,y)$ is $\mu$-primitive.
If  $h^\prime$ is primitive, then $\gamma (h^\prime)=1$.
\endproclaim

The crucial statement is

\proclaim{Theorem 3.1.3}
Let $X$ be a K3 surface with $\rho (X)=2$, and $H$ its polarization
of the degree $H^2=2rs$ with $\gamma (H)=1$.
Let $Y$ be the moduli space
of sheaves on $X$ with the isotropic Mukai vector $v=(r,H,s)$
and the canonical $nef$ element $h=(-a,0,b)\mod \bz v$.

As in Proposition 3.1.1, we consider the invariants
$\pm \mu \mod 2rs \in (\bz/2rs)^\ast$ and $d\in \bn$ of $(N(X), H)$.
Thus, we have
$$
\gamma(H)=1,\ \det N(X)=-d\ \text{where\ } \ \mu^2 \equiv d\mod 4rs.
\tag{3.1.8}
$$

With notations of Propositions {3.1.1}, then
all elements $\widetilde{h}=(xH+y\delta)/2rs \in N(X)$ with square
$\widetilde{h}^2=2ab$ satisfying Theorem {2.2.1} are in one to one
correspondence with integral solutions $(x,y)$ of the equation
$$
x^2-dy^2=4a^2b^2 c^2
\tag{3.1.9}
$$
which satisfy the conditions (i), (ii), (iii), (iv) below:

(i) $x\equiv \mu y\mod 2rs$;

(ii) $x\equiv \pm 2abc\mod d$;

(iii) $(x,y)$ belongs to one of $a$-series or $b$-series of solutions
defined below:

$a$-series: $b|(x,y)$ and $x-\mu y\equiv 0\mod 2rsa$ where
$\mu^2\equiv d\mod 4rsa$,

$b$-series: $a|(x,y)$ and $x-\mu y\equiv 0\mod 2rsb$ where
$\mu^2\equiv d\mod 4rsb$;

(iv) the pair $(x,y)$ is $\mu$-primitive: there does not exist
a prime $l$ such that $l|(x,y)$ and $x-\mu y\equiv 0\mod 2rs\,l$
(we always have $l^2|ab$ for such $l$), i. e.
$$
\text{g.c.d}\left(x,y,{x-\mu y\over 2rs}\right)=1.
\tag{3.1.10}
$$

In particular, by Theorem {2.2.1}, for a general $X$ with $\rho(X)=2$,
we have $Y\cong X$, if and only if the equation $x^2-dy^2=4a^2b^2c^2$
has an integral solution $(x,y)$ satisfying conditions (i)---(iv)
above. Moreover, a $nef$ primitive
element $\tilde{h}=(xH+y\delta)/2rs$ with
$\tilde{h}^2=2ab$ defines the pair $(X,\tilde{h})$ which is isomorphic to the
$(Y,h)$, if and only if $(x,y)$ satisfies the conditions (ii) and (iii)
(it satisfies conditions (i) and (iv) since it corresponds
to a primitive element of $N(X)$).
\endproclaim

\demo{Proof} We denote
$$
H^\ast={H\over 2rs},\ \delta^\ast={\delta\over 2rsd}
\tag{3.1.11}
$$
where $K(H)=\bz \delta=H^\perp$ in $N(X)$.
Since $\mu H^\ast+\delta/(2rs)\in N(X)$, then
$H^\ast+\mu^{-1}\delta/(2rs)\in N(X)$ where $\mu^{-1}$ is
defined $\mod 2rs$. It follows that
$$
u^\ast(H)={\mu^{-1}\over 2rs}\delta+K(H)= \mu^{-1}d\delta^\ast+\bz \delta.
\tag{3.1.12}
$$

Let
$$
\widetilde{h}={xH+y\delta\over 2rs} \in N(X)
\tag{3.1.13}
$$
satisfies conditions of Theorem 2.2.1. Then
$\widetilde{h}^2=2ab$. This is equivalent to
$x,\,y\in \bz$, $x\equiv \mu y\mod 2rs$ and $x^2-dy^2=4a^2b^2c^2$.
We get (i). The $\wth$ is primitive, if and only if $(x,y)$ is
$\mu$-primitive. We get (iv).

Consider $K(\wth)=\wth^\perp$ in $N(X)$.
If $k=\a H+\b \delta \in K(\wth)$, then
$(k,\widetilde{h})=\a x-\b yd=0$ and
$(\a,\b)=\lambda(yd,x)$, $\lambda\in \bq$.
We have
$$
(\lambda(ydH+x\delta))^2=\lambda^2(2rs y^2d^2-2rs dx^2)=
2rs \lambda^2d(y^2d-x^2)=-2rs \lambda^2d(4c^2a^2b^2).
$$
For $\lambda=1/2rs$ we get
$$
g=(ydH+x\delta)/2rs.
\tag{3.1.14}
$$
The element $g\in N(X)$ since $yd$, $x$
are integers and $yd-\mu x\equiv \mu(\mu y-x)\equiv 0\mod 2rs$.
We have $g^2=-2abd$. If $\wth$ is primitive and $\gamma (\wth)=1$,
then $K(\wth)=\bz g$. Then $N(X)=[\wth,g,(\nu\wth+g)/2ab]$ where
$\nu$ is a similar invariant for $\wth$ as $\mu$ for $H$. We denote
$$
\wth^\ast={\wth\over 2ab}\ \ \text{and}\ \ g^\ast={g\over 2abd}.
\tag{3.1.15}
$$
Then
$$
u^\ast(\wth)={\nu^{-1}g\over 2ab}=\nu^{-1}dg^\ast+K(h).
\tag{3.1.16}
$$
There exists a unique (up to $\pm 1$)
embedding $f:K(H)=\bz\delta\to K(\wth)=\bz g$
of one-dimensional lattices. It is given by $f(\delta)=\pm cg$.
Its dual is defined by
$f^\ast(g^\ast)=\pm c\delta^\ast$. We have $f^\ast(K(h))=
\bz f^\ast(g)=\bz (\delta/c)=\bz 2abcu^\ast(H)$ because of
\thetag{3.1.12}. This gives the first part of \thetag{2.2.13}.

Moreover,
$$
f^\ast(u^\ast (\wth))+f^\ast (K(h))=
\nu^{-1}df^\ast(g^\ast)+f^\ast(K(h))=
$$
$$
\pm \nu^{-1}cd\delta^\ast+f^\ast(K(h))=
\pm \nu^{-1}cu^\ast(H)+[K(H),2abcu^\ast(H)].
$$
Thus, the second part of \thetag{2.2.13} is
$\pm \nu^{-1}\equiv m(a,b)\mu^{-1}\mod 2ab$. Equivalently,
$$
\nu\equiv \pm m(a,b)\mu \mod 2ab.
\tag{3.1.17}
$$
Thus, for $\nu$ given by \thetag{3.1.17} one has
$$
\nu \widetilde{h}^\ast+g^\ast=\left(\left((\nu x+dy)H+
(\nu y+x)\delta\right)/2rs\right)/2ab \in N(X).
\tag{3.1.18}
$$
This is equivalent to
$$
\cases
\nu x+dy\equiv 0\mod 2ab\\
\nu y+x\equiv 0\mod 2ab\\
\nu x+dy\equiv \mu (\nu y+x)\mod 2rs2ab
\endcases .
\tag{3.1.19}
$$
For example, assume that $\nu\equiv m(a,b)\mu\mod 2ab$.
Then $\nu\equiv -\mu\mod 2a$ and $\nu\equiv \mu \mod 2b$.

From \thetag{3.1.19}, we then get $\mod 2a$ that
$-\mu x+y\mu^2\equiv 0\mod 2a$, $-\mu y+x\equiv 0\mod 2a$
and $2\mu x\equiv y(d+\mu^2)\mod 2rs 2a$. First two congruences
are valid since $x\equiv \mu y\mod 2rs$. The last one gives
$$
\mu {x-\mu y\over 2rs}\equiv {d-\mu^2\over 4rs}y\mod a \,.
\tag{3.1.20}
$$
We have $d\equiv \mu^2\mod 4rs$. Since $\mu\mod 2rs$ is invertible
and $a|r$, there exists the unique its lifting $\mu\mod 2rs a$
such that $d\equiv \mu^2\mod 4rs a$.
Then \thetag{3.1.20} is equivalent to  $x-\mu y\equiv 0\mod 2rsa$.

Now let us consider \thetag{3.1.19} $\mod 2b$.
First two congruences in \thetag{3.1.19} give
$\mu x+yd\equiv 0\mod 2b$, $\mu y+ x\equiv 0\mod 2b$.
Since $x\equiv \mu y\mod 2b$, they are equivalent to
$x\equiv y\equiv 0\mod b$. Third congruence in \thetag{3.1.19} gives
$yd\equiv \mu^2 y\mod 2rs2b$. It is valid since $b|y$ and
$\mu^2\equiv d\mod 4rs$. Thus, $(x,y)$ belongs to the $a$-series in (iii).
If $\nu\equiv -m(a,b)\mu\mod 2ab$, then similarly $(x,y)$
belongs to the $b$-series in (iii). Thus, condition \thetag{2.2.13} is
equivalent to (iii).

By \thetag{2.2.9},
$$
\split
&\kappa (H)(\delta^\ast)=-2rs(\delta^\ast\cdot
u^\ast(H))H^\ast+\delta^\ast+K(H)=\\
&-2rs(\delta/(2rsd)\cdot \mu^{-1}\delta/(2rs))H^\ast+\delta^\ast+K(H)=
\mu^{-1}H^\ast+\delta^\ast+K(H).
\endsplit
$$
Thus, we have
$$
\kappa (H)(\delta^\ast)=\mu^{-1}H^\ast+\delta^\ast+K(H).
\tag{3.1.21}
$$
Similarly,
$$
\kappa (\wth)(g^\ast)=\nu^{-1}g^\ast+K(\wth)=
\pm m(a,b)\mu^{-1}g^\ast+K(\wth).
\tag{3.1.22}
$$
Here we use \thetag{3.1.17}.

The condition \thetag{2.2.14} gives then
$$
\kappa (\wth)(g^\ast)=\pm c \kappa (H)(\delta^\ast).
\tag{3.1.23}
$$
From \thetag{3.1.21} and \thetag{3.1.22} we then get
$$
\split
&\kappa (\wth)(g^\ast)=
\pm m(a,b)\mu^{-1}{\wth\over 2ab}+ {g\over 2abd}+N(X)=\\
&\pm m(a,b)\mu^{-1}{xH+y\delta\over 2rs2ab}+
{ydH+x\delta\over 2rs2abd}+N(X)=\\
&{\pm m(a,b)\mu^{-1}x+y\over 2ab}H^\ast+
{\pm m(a,b)\mu^{-1}dy+x\over 2ab}\delta^\ast+N(X)\,.
\endsplit
\tag{3.1.24}
$$
We have $A_{N(X)}=N(X)^\ast/N(X)\cong (\bz/d\bz)$. The $\kappa(H)$
gives an epimorphism of $\bz\delta^\ast$ onto $A_{N(X)}$
with the kernel which is
$d\delta^\ast$. It follows that $\kappa (\wth)(g^\ast)=
\pm c\kappa (H)(\delta^\ast)$ is equivalent to
$$
{\pm m(a,b)\mu^{-1}dy+x\over 2ab}
\equiv \pm c\mod d.
\tag{3.1.25}
$$
Thus, $\pm m(a,b)\mu^{-1}yd+x\equiv \pm 2abc\mod 2ab d$.
Since $\text{g.c.d}(2abc,d)=1$, it is equivalent to
$x\equiv \pm 2abc \mod d$ and
$\pm m(a,b)\mu^{-1}yd+x\equiv 0 \mod 2ab$. The last congruence
is valid because of $d\equiv \mu^2\mod 2ab$, \thetag{3.1.17} and
\thetag{3.1.19}.

This finishes the proof.
\enddemo

Now let us consider the $a$-series of solutions of
\thetag{3.1.9} satisfying the main congruence (ii).
Since $b|(x,y)$ and $\text{g.c.d}(b,d)=1$,
we have $(x,y)=(bx_1,by_1)$ where $(x_1,y_1)$
satisfy
$$
x_1^2-dy_1^2=4a^2c^2\ \ \text{and}\ \ x_1\equiv \pm 2ac\mod d.
\tag{3.1.26}
$$
Let us find all solutions of \thetag{3.1.26}. We apply the main
trick of \cite{4} (see also \cite{5}).
Considering $\pm (x_1,\,y_1)$, we can assume that
$x_1\equiv 2ac \mod d$.
Then $x_1=2ac-kd$ where $k\in \bz$. We have
$4a^2c^2-4ackd+k^2d^2-dy_1^2=4a^2c^2$. Thus,
$$
d={y_1^2+4ack \over k^2}.
\tag{3.1.27}
$$
Let $l$ be prime. Like in \cite{4}, it is easy to see that if
$l^{2t+1}\vert k$ and $l^{2t+2}\not| k$, then $l\vert 4ac$.
It follows that $k=-\alpha q^2$ where $\alpha|2ac$ and $\alpha$ is
square-free, $q\in \bz$. It follows that
\thetag{3.1.27} is equivalent to
$$
\alpha | 2ac \ \text{is square-free,}\ \ p^2-dq^2={4ac \over \alpha}.\
\tag{3.1.28}
$$
Then $k=-\alpha q^2$, $y=\alpha pq$. Here $\alpha$ can be negative.
Thus, (integral) solutions
$\alpha$, $(p,q)$ of \thetag{3.1.28} give all solutions
$$
(x_1,y_1)=\pm (2ac+\alpha dq^2,\ \alpha pq)
\tag{3.1.29}
$$
of \thetag{3.1.26}.
We call them {\it associated solutions.} Thus,
all solutions $(x_1,y_1)$ of \thetag{3.1.26} are associated solutions
\thetag{3.1.29} to all solutions $\alpha$, $(p,q)$ of \thetag{3.1.28}.
If one additionally assumes that $q\ge 0$, then $(x,y)$ and
$\alpha$, $(p,q)$ are in one to one correspondence, by our construction.

Now let us consider associated solutions \thetag{3.1.29} which satisfy
the additional condition $bx_1\equiv \mu by_1\mod 2rs$. (It is the
condition (i) of Theorem 3.1.3.) Equivalently,
$x_1\equiv \mu y_1\mod 2ac^2$. By \thetag{3.1.29}, this is equivalent
$$
2ac+\alpha dq^2\equiv \mu\alpha pq\mod 2ac^2.
\tag{3.1.30}
$$
We also have
$$
4ac+\alpha dq^2=\alpha p^2.
\tag{3.1.31}
$$
Taking \thetag{3.1.31} - \thetag{3.1.30}, we get that \thetag{3.1.30} is
equivalent to
$$
2ac\equiv \alpha p(p-\mu q)\mod 2ac^2.
\tag{3.1.32}
$$
Taking 2\thetag{3.1.30} - \thetag{3.1.31}, we get that \thetag{3.1.30}
is equivalent to
$$
\alpha p^2+\alpha dq^2-2\alpha \mu pq\equiv 0\mod 4ac^2.
\tag{3.1.33}
$$
Since $d\equiv \mu^2\mod 4abc^2$, this is equivalent to
$$
\alpha (p-\mu q)^2\equiv 0\mod 4ac^2.
\tag{3.1.34}
$$
Since $\alpha$ is square-free, it follows that
$$
2c\,\vert p-\mu q.
\tag{3.1.35}
$$
From \thetag{3.1.32}, we then get
$$
2a\equiv \alpha p \left({p-\mu q\over c}\right)\mod {2ac}
\tag{3.1.36}
$$
where $(p-\mu q)/c$ is an integer. Since $\alpha |2ac$,
it follows that $\alpha | 2a$. Thus, we get
$$
\alpha \vert 2a,\ \ c\vert (p-\mu q)\ \ \text{and}\ \
{2a\over \alpha}\equiv p \left({p-\mu q\over c}\right)\mod
{{2ac\over \alpha}}.
\tag{3.1.37}
$$
Thus, we have proved (see \thetag{3.1.34}) that
the condition $x_1\equiv \mu y_1\mod 2ac^2$ is equivalent to
$$
\alpha \vert 2a\ \ \text{and}\ \
(p-\mu q)^2\equiv 0\mod {{4ac^2\over \alpha}}.
\tag{3.1.38}
$$
Let us prove that actually $\alpha \vert a$. Otherwise, the
equation $p^2-dq^2=4ac/\alpha$ has no solutions. Really,
assume that $2a/\alpha$ is odd. By \thetag{3.1.38} we get
$(p-\mu q)^2 \equiv 0 \mod 2 c^2$. It follows
$p-\mu q \equiv 0 \mod 2c$ and $p+\mu q\equiv 0\mod 2$. Then
$p^2-\mu^2 q^2\equiv 0 \mod 4 c$. Since $\mu^2\equiv d\mod 4ac^2$,
we then get $p^2-d q^2 \equiv 0 \mod 4c$. On the other hand, if
$2a/\alpha$ is odd, then $p^2-dq^2=4ac/\alpha \equiv 2c\mod 4c$.
We get a contradiction.

Now we add the condition (iv) and the last condition
of the $a$-series $x-\mu y\equiv 0\mod 2rsa$ of Theorem 3.1.3 (here
we assume that $\mu^2\equiv d\mod 4rsa$). For
$(x,y)=(bx_1,by_1)$
we then get
$$
x_1-\mu y_1\equiv 0 \mod 2a^2c^2\ \text{and}\
\text{g.c.d}(bx_1,\,by_1,\,{x_1-\mu y_1\over 2ac^2})=1.
\tag{3.1.39}
$$
By \thetag{3.1.29}, we then have $\alpha|(x_1,y_1,(x_1-\mu y_1)/2ac^2)$.
It follows the fundamental result:
$$
\alpha=\pm 1 .
\tag{3.1.40}
$$
The condition $x_1-\mu y_1\equiv 0 \mod 2a^2c^2$ gives
$$
2ac+\alpha dq^2\equiv \mu\alpha pq\mod 2a^2c^2.
\tag{3.1.41}
$$
As always
$$
-4ac+\alpha p^2-\alpha dq^2=0.
\tag{3.1.42}
$$
Taking 2\thetag{3.1.41}+\thetag{3.1.42}, we get that \thetag{3.1.41}
is equivalent to
$$
\alpha dq^2+\alpha p^2\equiv 2\mu\alpha pq\mod 4a^2c^2.
\tag{3.1.43}
$$
Since $\alpha=\pm 1$, this is equivalent to
$$
dq^2+p^2\equiv 2\mu pq\mod 4a^2c^2.
\tag{3.1.44}
$$
Since $\mu^2\equiv d\mod 4a^2c^2$, the \thetag{3.1.44} is equivalent to
$$
(p-\mu q)^2\equiv 0\mod 4a^2c^2
\tag{3.1.45}
$$
which is equivalent to
$$
p-\mu q\equiv 0\mod 2ac.
\tag{3.1.46}
$$
Thus, the congruence $x_1-\mu y_1\equiv 0 \mod 2a^2c^2$ in \thetag{3.1.39}
is equivalent to \thetag{3.1.46}. Now let us analyze the
condition of primitivity
$$
\text{g.c.d}(bx_1,\,by_1,\,{x_1-\mu y_1\over 2ac^2})=1
\tag{3.1.47}
$$
in \thetag{3.1.39}. The condition \thetag{3.1.47} is
equivalent to the fact that there does not exist a prime $l$ such that
$$
l|(bx_1,\,by_1,\,(x_1-\mu y_1)/2ac^2)\ \text{where always}\ l^2|ab.
\tag{3.1.48}
$$

Assume $l^2|a$. Since $a|(x_1-\mu y_1)/2ac^2$, \thetag{3.1.48}
is equivalent to $l|(x_1,\,y_1)$. We have $(x_1,\,y_1)=
\pm (2ac+\alpha dq^2,\, \alpha pq)=
\pm (-2ac+\alpha p^2,\,\alpha pq)$. It follows $l|(p,q)$. Thus,
the $l^2|a$ in \thetag{3.1.48} does not exist, if and only if
$l|(p,q)$ is impossible for $l^2|a$.

Assume $l^2|b$. Then $l|(bx_1,\,by_1,\,(x_1-\mu y_1)/2ac^2)$ is
equivalent to $x_1-\mu y_1\equiv 0\mod 2ac^2 l$. This is
equivalent to
$x_1-\mu y_1\equiv 0\mod 2a^2c^2 l$ because we know that
$x_1-\mu y_1\equiv 0\mod 2a^2c^2$ and $\text{g.c.d}(l,a)=1$. From
\thetag{3.1.29}, we get
$$
2ac+\alpha d q^2\equiv \mu \alpha pq\mod 2a^2c^2 l .
\tag{3.1.49}
$$
By $-4ac+\alpha p^2-\alpha d q^2=0$ and $\alpha=\pm 1$,
we then get like above that \thetag{3.1.49} is equivalent to
$$
dq^2+p^2-2\mu pq \equiv 0 \mod 4a^2c^2 l.
\tag{3.1.50}
$$
Since $\mu^2\equiv d\mod 4abc^2 a$ and $l|b$,
the \thetag{3.1.50} is equivalent to
$$
(p-\mu q)^2\equiv 0\mod 4a^2c^2\,l.
\tag{3.1.51}
$$
The \thetag{3.1.51} is equivalent to
$$
p-\mu q\equiv 0\mod 2ac\,l.
\tag{3.1.52}
$$
Thus, we finally get that \thetag{3.1.39} is equivalent to
$$
\cases
p\equiv \mu q\mod 2ac,\\
\text{g.c.d}(a,p,q)=1,\\
p\not\equiv \mu q\mod 2acl\ \text{for any prime}\ l|b.
\endcases
\tag{3.1.53}
$$

We get the final result:

\proclaim{Theorem 3.1.4}
With conditions of Theorem 3.1.3,
for a general K3-surface $X$ with $\rho (X)=2$ and $\gamma(H)=1$,
we have $Y\cong X$, if and only if for one of $\alpha=\pm 1$ there exists
an integral solution $(p,q)$ for one of $a$-series or $b$-series:

$a$-series:
$$
\cases
p^2-dq^2=4ac/\alpha \\
p\equiv \mu q \mod 2ac\\
\text{g.c.d}(a,p,q)=1\\
p\not\equiv \mu q\mod 2acl\ \text{for any prime}\  l|b
\endcases .
\tag{3.1.54}
$$
A solution $(p,q)$ of \thetag{3.1.54} gives a solution $(x,y)$
of Theorem 3.1.3 as associated solution
$$
(x,y)=\pm (2abc+\alpha bdq^2,\ \alpha  bpq).
\tag{3.1.55}
$$

For $b$-series, one should just change $a$ and $b$ places:

$b$-series:
$$
\cases
p^2-dq^2=4bc/\alpha \\
p\equiv \mu q \mod 2bc\\
\text{g.c.d}(b,p,q)=1\\
p\not\equiv \mu q\mod 2bcl\ \text{for any prime}\  l|a
\endcases .
\tag{3.1.56}
$$
A solution $(p,q)$ of \thetag{3.1.56} gives a solution $(x,y)$
of Theorem 3.1.3 as associated solution
$$
(x,y)=\pm (2abc+\alpha adq^2,\ \alpha apq).
\tag{3.1.57}
$$
All $\wth$ (or solutions $(x,y)$) of Theorem 3.1.3 are given by
\thetag{3.1.55} and \thetag{3.1.57} as associated solutions.
\endproclaim

Now let us interpret solutions $(p,q)$ of Theorem 3.1.4 as
elements of the Picard lattice $N(X)$. Further $\pm 1$ denotes
$\alpha$. If there is an ambiguity with other $\pm$, we introduce
$\alpha$ (e.g. see \thetag{3.1.63} below).

Let us consider a solution $(p,q)$ of the $a$-series \thetag{3.1.54}.
The congruence $p\equiv \mu q\mod 2ac$ is equivalent to
$$
tp\equiv t \mu q\mod 2abc^2\, ,
\tag{3.1.58}
$$
where
$$
t=bc .
\tag{3.1.59}
$$
The congruence \thetag{3.1.58} is equivalent to
$$
h_1=\frac{t(pH+q\delta)}{2abc^2} \in N(X)\  \text{and}\
h_1\cdot H\equiv 0\mod t.
\tag{3.1.60}
$$
We have $p^2-dq^2=\pm 4ac$ is
equivalent to $(h_1)^2=t^2(p^2-dq^2)/(2abc^2)=\pm 2bc$.
Moreover, $H\cdot h_1\equiv 0\mod bc$.

The condition $\text{g.c.d}(a,p,q)=1$ can be interpreted as
$l$ does not divide $p=H\cdot h_1/t$, if $l^2|a$. Really, if
$l|a$ and $l|p$, then $l|q$ since $p^2-dq^2=\pm 4ac$. Thus, we
get $l|(H\cdot h_1)/bc$ is impossible for $l^2|a$.

Assume that $p\equiv \mu q\mod 2acl$ for a prime $l|b$. This
is equivalent to $pt/l\equiv \mu qt/l\mod 2abc^2$. Equivalently,
$h_1/l=\left((pt/l)H+(qt/l)\delta\right)/(2abc^2)\in N(X)$. Thus,
such a prime $l$ does not exist, if and only if $h_1$ is not divisible
by $l|b$ in $N(X)$. It is enough to consider $l^2|b$.

Changing $a$ and $b$ places, we get similar statements for the $b$-series.

Thus, we finally get

\proclaim{Theorem 3.1.5}
With notations and conditions of Theorem 3.1.3, for a general K3-surface
$X$ with $\rho (X)=2$ and $\gamma(H)=1$ for $H\in N(X)$,
we have $Y\cong X$, if and only if
at least for one of signs $\pm$ (i. e. $\alpha=\pm 1$)
there exists $h_1\in N(X)$ such that $h_1$ belongs to
one of $a$-series or $b$-series:

$a$-series:
$$
h_1^2=\pm 2bc,\  H\cdot h_1\equiv 0\mod bc,\
H\cdot h_1\not\equiv 0\mod bcl_1,\
h_1/l_2\not\in  N(X)
\tag{3.1.61}
$$
where $l_1$ and $l_2$ are any primes such that $l_1^2|a$ and $l_2^2|b$;

$b$-series:
$$
h_1^2=\pm 2ac,\  H\cdot h_1\equiv 0\mod ac,\
H\cdot h_1\not\equiv 0\mod acl_1,\
h_1/l_2\not\in  N(X)
\tag{3.1.62}
$$
where $l_1$ and $l_2$ are any primes such that $l_1^2|b$ and $l_2^2|a$.
\endproclaim

\remark{Important Remark} By Theorem 3.1.3 and the formulae \thetag{3.1.55},
\thetag{3.1.57}
for the associated solution, we get that
$$
\pm \widetilde{h}=
\cases
&-\frac{H}{c} +\frac{\alpha (H\cdot h_1)h_1}{bc^2},\ \text{if $h_1$
is from $a$-series,}\\
&-\frac{H}{c} + \frac{\alpha (H\cdot h_1)h_1}{ac^2},\ \text{if $h_1$
is from $b$-series}
\endcases
\tag{3.1.63}
$$
belongs to $N(X)$ and
$$
(Y,h)\cong \left(X,\pm w(\widetilde{h})\right)\ \text{for some}\
w\in W^{(-2)}(N(X)).
\tag{3.1.64}
$$
It describes moduli $(Y,h)$ of sheaves in terms of $X$ where
$h=(-a,0,b)\mod \bz v$.
\endremark

Applying additionally Theorem 2.2.1, we get the following simple
sufficient condition when $Y\cong X$ which is valid for $X$ with
any $\rho (X)$. This is one of the main results of the paper.

\proclaim{Theorem 3.1.6}
Let $X$ be a K3 surface with a polarization $H$ such that
$H^2=2rs$, $r,\,s\ge 1$, the Mukai vector $(r,H,s)$ is primitive
and $Y$ the moduli of sheaves on $X$ with the Mukai vector
$v=(r,H,s)$.

Then $Y\cong X$, if at least for one of signs $\pm$ (i. e. $\alpha=\pm 1$)
there exists $h_1\in N(X)$ such that elements $H$, $h_1$ are contained in a
2-dimensional sublattice $N\subset N(X)$ with $H\cdot N=\bz$ and $h_1$
belongs to one of $a$-series or $b$-series described below:

$a$-series:
$$
h_1^2=\pm 2bc,\  H\cdot h_1\equiv 0\mod bc,\
H\cdot h_1\not\equiv 0\mod bcl_1,\
h_1/l_2\not\in  N(X)
\tag{3.1.65}
$$
where $l_1$ and $l_2$ are any primes such that $l_1^2|a$ and $l_2^2|b$;

$b$-series:
$$
h_1^2=\pm 2ac,\  H\cdot h_1\equiv 0\mod ac,\
H\cdot h_1\not\equiv 0\mod acl_1,\
h_1/l_2\not\in  N(X)
\tag{3.1.66}
$$
where $l_1$ and $l_2$ are any primes such that $l_1^2|b$ and $l_2^2|a$.

We also have the formulae \thetag{3.1.63} and \thetag{3.1.64} for
$(Y,h)$ where $h=(-a,0,b)\mod \bz v$.

The conditions above are necessary for $\gamma(H)=1$ and
$Y\cong X$, if $\rho(X)\le 2$ and $X$ is a general K3 surface with
its Picard lattice, i. e. the automorphism group of the transcendental
periods $(T(X), H^{2,0}(X))$ of $X$ is $\pm 1$.
\endproclaim

\head
4. Divisorial conditions on moduli of $(X,H)$ which imply
$\gamma (H)=1$ and $Y\cong X$
\endhead

Further we use the following notations. We fix:

$r,s\in \bn$, $r,s\ge 1$, $r=ca, s=cb$ where $c=\text{g.c.d.}(r,s)$
and then $\text{g.c.d.}(a,b)=1$;

$\overline{\mu}=\{\mu,-\mu\}\subset (\bz/ 2rs)^\ast$.

For any choice of the sign $\pm$ (equivalently, $\alpha=\pm 1$),
we denote
$$
\align
\Da(r,s;a)^{\overline{\mu}}_\pm
=\{
& d\in \bn \ | \
d\equiv \mu^2 \mod 4abc^2, \
\exists \
(p,q) \in \bz \times \bz:\  p^2-dq^2=\pm 4ac,\  \\
&p-\mu q \equiv 0 \mod 2ac,\
\text{g.c.d}(a,p,q)=1,\\
&p-\mu q\not\equiv 0 \mod 2ac l\
\text{for any prime}\ l\ \text{such that}\ l^2|b\}.
\tag{4.1}
\endalign
$$
The sets $\Da(r,s;a)^{\overline{\mu}}_\pm$ are called the $a$-series.

For any choice of the sign $\pm$ we denote
$$
\align
\Da(r,s;b)^{\overline{\mu}}_\pm
=\{
& d\in \bn \ | \
d\equiv \mu^2 \mod 4abc^2, \
\exists \
(p,q) \in \bz \times \bz:\  p^2-dq^2=\pm 4bc,\  \\
&p-\mu q \equiv 0 \mod 2bc,\
\text{g.c.d}(b,p,q)=1,\\
&p-\mu q\not\equiv 0 \mod 2bc l\
\text{for any prime}\ l\ \text{such that}\ l^2|a\}.
\tag{4.2}
\endalign
$$
The sets $\Da(r,s;b)^{\overline{\mu}}_\pm$ are called the $b$-series.

Let $X$ be a K3 surface with a polarization $H$
of the degree $H^2=2rs$ and a primitive Mukai vector $(r,H,s)$.
The moduli space $Y$ of sheaves
on $X$ with the Mukai vector $v=(r,H,s)$ has the canonical $nef$
element
$h=(-a,0,b) \mod \bz v$ with $h^2=2ab$.

If $\gamma (H)=1$, then $\rho (X)\ge 2$.
Since the dimension of the
moduli space of $(X,H)$ is equal to $20-\rho (X)$, it follows that describing
general $(X,H)$ with $\rho (X)=2$, $\gamma (H)=1$ and $Y\cong X$,
we at the same time describe all possible divisorial conditions
on moduli of $(X,H)$ which imply  $\gamma (H)=1$ and
$Y\cong X$.
They are described by the invariants of the pairs $(N(X),H)$
where $\rk N(X)=2$. Here we use results of \cite{10,11}
(see also \cite{3}) about irreducibility of moduli of K3 surfaces
with a fixed Picard lattice.

Using Theorem 3.1.4, we get

\proclaim{Theorem 4.1}
With the above notations all possible divisorial conditions on
moduli of polarized K3 surfaces $(X,H)$ with a polarization
$H$ of the degree $H^2=2rs$, which imply $\gamma (H)=1$ and
$Y\cong X$ are labelled by the set
$$
{\Cal Div}(r,s)={\Cal Div}(r,s;a)\bigcup {\Cal Div}(r,s;b),
\tag{4.3}
$$
where
$$
\split
{\Cal Div}(r,s;a)=\{(d,\,\overline{\mu})\ |\
\overline{\mu}=\{\mu, -\mu\}\subset
(\bz/2rs)^\ast,\ d\in \Da(r,s;a)^{\overline{\mu}}\},\\
{\Cal Div}(r,s;b)=\{(d,\,\overline{\mu})\ |\
\overline{\mu}=\{\mu, -\mu\}\subset
(\bz/2rs)^\ast,\ d\in \Da(r,s;b)^{\overline{\mu}}\}
\endsplit
\tag{4.4}
$$
and
$$
\Da(r,s;a)^{\overline{\mu}}=\Da(r,s;a)^{\overline{\mu}}_+
\bigcup \Da(r,s;a)^{\overline{\mu}}_-,\ \
\Da(r,s;b)^{\overline{\mu}}=\Da(r,s;b)^{\overline{\mu}}_+
\bigcup \Da(r,s;b)^{\overline{\mu}}_-.
\tag{4.5}
$$

The sets
$$
\Da(r,s;a)^{ \overline{\mu}}_{\pm}\ \text{and}\
\Da(r,s;b)^{ \overline{\mu}}_{\pm}
\tag{4.6}
$$
are infinite if they are not empty.

The set
$\Da(r,s;a)^{ \overline{\mu}}_{\pm}$
is not empty, if and only if the congruence
$$
\mu q t +act^2\equiv \pm 1 \mod bcq^2
\tag{4.7}
$$
has a solution $q\in\bn$, $t\in \bz$. E.g. it is true for $q=1$,
if $b|c^2$. The set $\Da(r,s;b)^{ \overline{\mu}}_{\pm}$
is not empty, if and only if
the congruence
$$
\mu q t +bct^2\equiv \pm 1 \mod acq^2
\tag{4.8}
$$
has a solution for $q\in \bn$, $t\in \bz$.
E. g. it is true for $q=1$, if $a|c^2$.
See other general results of this type in Theorems 4.5 and 4.6 below.

The set ${\Cal Div}(r,s)$ is always non-empty and then infinite: i. e.
at least one of sets $\Da(r,s;a)^{ \overline{\mu}}_{\pm}$ or
$\Da(r,s;b)^{ \overline{\mu}}_{\pm}$ is not empty for one of
$\overline{\mu}$ and $\pm$.
Its subset ${\Cal Div}(r,s;a)\subset{\Cal Div}(r,s)$
is non-empty and then infinite, if
either $abc$ is odd or $ac$ is even; its subset
${\Cal Div}(r,s;b)\subset{\Cal Div}(r,s)$ is non-empty
and then infinite, if either $abc$ is odd or $bc$ is even.
\endproclaim

Below we prove statements of Theorem 4.1 which don't follow
from Theorem 3.1.4 straightforward.
These considerations are important as itself, and they give more
exact and strong statements as we have formulated in Theorem 4.1.

Further we consider the $a$-series. The same results are valid
for the $b$-series changing $a$ and $b$ places.
For $q\in \bn$ such that $\text{g.c.d}(a,q)=1$ we introduce
$$
\align
\Da(r,s;a)^{\overline{\mu}}_\pm(q)
=\{
& d\in \bn \ | \
d\equiv \mu^2 \mod 4abc^2, \
\exists \
p \in \bz:\  p^2-dq^2=\pm 4ac,\  \\
&p-\mu q \equiv 0 \mod 2ac,\
\text{g.c.d}(a,p,q)=1,\\
&p-\mu q\not\equiv 0 \mod 2ac l\
\text{for any prime}\ l\ \text{such that}\ l^2|b\}.
\tag{4.9}
\endalign
$$
By definition,
the set $\Da(r,s;a)^{\overline{\mu}}_\pm(q)=\emptyset$, if
$\text{g.c.d}(a,q)>1$.

Obviously,
$$
\Da(r,s;a)^{\overline{\mu}}_\pm=\bigcup_{q\in \bn}
{\Da(r,s;a)^{\overline{\mu}}_\pm (q)}.
\tag{4.10}
$$
We have

\proclaim{Lemma 4.2} For any $\pm$ and any $q\in \bn$ the set
$\Da(r,s;a)^{\overline{\mu}}_\pm (q)$ is either empty or
infinite. It is infinite, if and only if the system of congruences
$$
\cases
p\equiv \mu q\mod 2ac,\\
p^2\equiv \pm 4ac+\mu^2q^2\mod {4abc^2q^2},\\
p\not\equiv \mu q\mod 2acl
\endcases
\tag{4.11}
$$
has a solution. Here $l$ is any prime such that $l^2|b$,
the condition $p\not\equiv \mu q\mod 2acl$ must satisfy for
all these primes. The set $\Da(r,s;a)^{\overline{\mu}}_\alpha(q)$
consists of all $d=(p^2\mp 4ac)/q^2$ where $p$ satisfies the above system
of congruences and $p^2\mp 4ac>0$.
\endproclaim

\demo{Proof} It at once follows from definitions using that
$d=(p^2\mp4ac)/q^2>0$. Clearly, if $p=p_0$ satisfies the
system \thetag{4.11}, then $p=p_0+t(4abc^2q^2)$, $t\in \bz$,
satisfies \thetag{4.11}, and $p^2\mp 4ac >0$ for $|t|>>0$.
It gives an infinite subset of
$\Da(r,s;a)^{\overline{\mu}}_\pm (q)$. It finishes the proof.
\enddemo

We denote
$$
\alpha=\pm 1.
\tag{4.12}
$$
Let us rewrite \thetag{4.11}. We have
$p=\mu q+2act$, $t\in \bz$ where $\text{g.c.d}(t,b)=1$.
Then $(\mu q+2act)^2
\equiv 4\alpha ac+\mu^2q^2\mod~{4abc^2q^2}$.
It follows $\mu q t +act^2\equiv \pm 1 \mod bcq^2$. Then
$\text{g.c.d}(t,b)=1$ and $\text{g.c.d}(q,a)=1$. It follows
that \thetag{4.11} is equivalent to the congruence
$$
\mu q t +act^2\equiv \pm 1 \mod bcq^2.
\tag{4.13}
$$
A solution $t\in \bz$ of \thetag{4.13} gives the infinite
series of elements from
${\Da(r,s;a)^{\overline{\mu}}_\pm (q)}$ which is
$$
{\left(\mu q +2ac(t+bc q^2 m)\right)^2\mp 4ac
\over q^2}\ >\ 0,\ \ m\in \bz\,.
\tag{4.14}
$$
Thus, additionally to Lemma 4.2 we get

\proclaim{Lemma 4.3} For $q\in \bn$ the set
$\Da(r,s;a)^{\overline{\mu}}_\pm(q)$ is non-empty, if and only if
the congruence \thetag{4.13} has a solution $t\in \bz$. Each such
solution $t$ generates an infinite series of elements from
$\Da(r,s;a)^{\overline{\mu}}_\pm (q)$ given by \thetag{4.14}. Together
they give all elements of $\Da(r,s;a)^{\overline{\mu}}_\pm(q)$.

In particular, the set $\Da(r,s;a)^{\overline{\mu}}_\pm$
is non-empty (and then it is infinite),
if and only if the congruence \thetag{4.13} has a solution
$(q,t)$ for $q\in \bn$ and $t\in \bz$.
\endproclaim

For $q=1$, the congruence \thetag{4.13} gives
$$
\mu t +act^2\equiv \alpha \mod~{bc}.
\tag{4.15}
$$
We have $\mu t\equiv \alpha \mod c$ and $t=\mu^{-1}\alpha+kc$, $k\in \bz$.
Simple calculations show that \thetag{4.15} is equivalent
to the quadratic equation
$$
k^2ac^2+k(\mu+2\alpha ac\mu^{-1})+a\mu^{-2}\equiv 0\mod b.
\tag{4.16}
$$
For $A\not=0$, the quadratic equation
$$
Ax^2+Bx+C\equiv 0\mod b
\tag{4.17}
$$
is equivalent to
$$
(2Ax+B)^2\equiv B^2-4AC\mod 4Ab .
\tag{4.18}
$$
It follows that \thetag{4.15} is equivalent to
$$
(2ac^2k+\mu+2\alpha ac\mu^{-1})^2\equiv (\mu^2+4\alpha ac)\mod 4abc^2.
\tag{4.19}
$$
In particular, \thetag{4.16} (equivalently \thetag{4.19}) has
a solution, if $c^2\equiv 0\mod b$. Then $A=ac^2\equiv 0\mod b$ and
$B=\mu+2\alpha ac\mu^{-1}$ is invertible $\mod b$. Thus, we get

\proclaim{Lemma 4.4} The set $\Da(r,s;a)^{\overline{\mu}}_{\pm}(1)$
is non-empty and then infinite, if and only if
the quadratic equation \thetag{4.16} (equivalently, \thetag{4.19})
has a solution $k$.

In particular, this is valid, if $b|c^2$.
\endproclaim

Using Lemma 4.4, we also get the following general result:

\proclaim{Theorem 4.5} The set
$\Da(r,s;a)^{\overline{\mu}}_{\alpha}(1)$
is non-empty for $\mu\equiv \theta^{-1}ac-\theta \alpha\mod 2abc^2$,
if $\theta\in (\bz/2abc^2)^\ast$ and
$\text{g.c.d}(2b,\theta^{-1}ac-\theta \alpha)=1$.

It follows that $\Da(r,s;a)^{\overline{\mu}}_{\alpha}(1)$
is non-empty at least for one of $\alpha$ and ${\overline{\mu}}$, if
$ac$ is even. In particular, the set ${\Cal Div}(r,s;a)$ is not empty,
if $ac$ is even.
\endproclaim

\demo{Proof} By conditions on $\theta$, we have that
$\mu\equiv \theta^{-1}ac-\theta\alpha\mod 2abc^2\in
(\bz/2abc^2)^\ast$. Then $\mu^2+4\alpha ac\equiv \nu^2\mod 4abc^2$
where $\nu\equiv \theta^{-1}ac+\theta \alpha\mod 2abc^2$.
Thus, \thetag{4.19} satisfies, if
$$
2ac^2k+\mu+2\alpha\,ac\mu^{-1}\equiv -\nu \mod 2abc^2
$$
has a solution $k$. Thus, we should have
$$
2ac^2k+(\theta^{-1}ac-\theta\alpha)+{2\alpha ac\over
\theta^{-1}ac-\theta \alpha}\equiv
-(\theta^{-1}ac+\theta \alpha)\mod 2abc^2.
$$
This gives
$$
(\theta^{-1}ac-\theta\alpha)2ac^2k\equiv -2\theta^{-2}a^2c^2\mod 2abc^2.
$$
Equivalently,
$$
(\theta^{-1}ac-\theta\alpha)k\equiv -\theta^{-2}a\mod b.
$$
It has a solution $k$ since
$\text{g.c.d}(\theta^{-1}ac-\theta\alpha, b)=1$.
By Lemma 4.4, the set \newline 
$\Da(r,s;a)^{\overline{\mu}}_{\alpha}(1)$ is not then empty.

Assume that $ac$ is even. We can find $\alpha$ such that for
any prime $p_i\vert 2b$ there exists
$x_i\mod p_i$ such that $x_i\not\equiv 0\mod p_i$ and
$-\alpha x_i^2+ac\not\equiv 0\mod p_i$.
Really, if $p_i=2$, we can take $x_i\equiv 1\mod 2$.
If $p_i=3$, we can take $x_i\equiv \pm 1\mod 3$.
If $ac\equiv 0\mod 3$, we can take any $\alpha=\pm 1$.  If
$ac\equiv 1\mod 3$, we take $\alpha=-1$. If $ac\equiv-1\mod 3$, we
take $\alpha=1$. If $p_i\ge 5$, the equation $-\alpha x^2+ac\equiv 0\mod p_i$
has not more then two solutions $\mod p_i$, and there always exists
$x_i\mod p_i$ with the above conditions.

By standard results, there exists
$x\mod 2b$ such that $x\equiv x_i\mod p_i$. It then follows that
$x\in (\bz/2b)^\ast$ and $(-\alpha x^2+ac,2b)=1$. By
standard results and considerations,  there exists
$\theta\in (\bz/2abc^2)^\ast$ such that
$\theta\equiv x\mod 2b$. Then $(-\alpha \theta^2+ac,2b)=1$ and
also $(-\alpha\theta+ac\theta^{-1},2b)=1$. For the found $\alpha$ and
$\theta$, the set $\Da(r,s;a)^{\overline{\mu}}_{\alpha}(1)$ is not empty,
if $\mu\equiv-\alpha\theta+ac\theta^{-1}\mod 2abc^2$.

This finishes the proof.
\enddemo

Considering Theorem 4.5 also for the $b$-series, we get that the
 set ${\Cal Div}(r,s)={\Cal Div}(r,s;a)\bigcup {\Cal Div}(r,s;b)$
is not empty, if at least one of $a$, $b$ or $c$ is even.

Thus, to prove that the set
${\Cal Div}(r,s)$ is not empty, we need to consider only the case when
$abc$ is odd. Then we have

\proclaim{Theorem 4.6} Assume that $abc$ is odd. Then the set
$\Da(r,s;a)^{\overline{\mu}}_{\alpha}(1)$
is non-empty for $\mu\mod 2abc^2\in (\bz/2abc^2)^\ast$, if
$\mu\equiv \theta^{-1}ac-\theta \alpha\mod abc^2$
where $\theta\in (\bz/abc^2)^\ast$ and
$\text{g.c.d}(b,\theta^{-1}ac-\theta \alpha)=1$.

It follows that $\Da(r,s;a)^{\overline{\mu}}_{\alpha}(1)$
is not empty at least for one of $\alpha$ and ${\overline{\mu}}$.
In particular, the set ${\Cal Div}(r,s;a)$ is not empty.
\endproclaim

\demo{Proof} The congruence \thetag{4.19} is the identity $\mod 4$.
Then, since $abc$ is odd, \thetag{4.19} is equivalent to
$$
(2ac^2k+\mu+2\alpha ac\mu^{-1})^2\equiv (\mu^2+4\alpha ac)\mod abc^2.
\tag{4.20}
$$
Thus, it is enough to consider $\mu \mod abc^2$. Then the same
considerations as for Theorem 4.5 prove the statement.
\enddemo

We can summarize the main results:

\proclaim{Theorem 4.7} For any $r$, $s$ the set
${\Cal Div}(r,s)$ of divisorial conditions on moduli is non-empty and then
infinite.

Its subset ${\Cal Div}(r,s;a)$ is non-empty and then
infinite, if either $abc$ is odd or $ac$ is even;
its subset ${\Cal Div}(r,s;b)$ is
non-empty and then infinite, if either $abc$ is odd or $bc$ is even.
\endproclaim

\head
5. The results from the viewpoint of $Y$.
\endhead

It is interesting and important to consider the same problems from
the viewpoint of $Y$. In this case, we fix $a,\,b\in \bn$ with
$\text{g.c.d}(a,b)=1$, and we fix $c\in \bn$. Further we denote
$r=ac$ and $s=bc$.

We consider a K3 surface $Y$ with a $nef$ element $h\in N(Y)$
such that $h^2=2ab$. We ask: {\it When $(Y,h)$ is the moduli space $\M$
of sheaves on $Y$ with the Mukai vector $v=(r,H,s)$ and the
canonical $nef$ element $\wth=(-a,0,b)\mod \bz v$ where
$H\in N(Y)$ is $nef$, $H^2=2rs$ and $\gamma (H)=1$?
More exactly, when there exists an isomorphism $\xi:\M\to Y$ such that
$\xi(\wth)=h$?}

This question is very similar to the previous ones for $X$, and,
in principle, the corresponding results can be deduced from the results
for $X$ above. But it is simpler to reconsider the previous calculations
from the viewpoint of $Y$. They are very similar but sometimes
very non-trivial. It is surprising that the results for $Y$ below are
so similar to the results for $X$.

Theorem 5.1 below is a reformulation of Theorem 2.2.1. We use
notations from Sect. 2 as in Theorem 2.2.1.

\proclaim{Theorem 5.1} Let $Y$ be a K3 surface with a $nef$ element
$h$ of the degree $h^2=2ab$.

Then $(Y,h)$ is (isomorphic to) the moduli space $(\M,\wth)$ of
coherent sheaves on $Y$ with the Mukai vector
$v=(r,H,s)$, where $H^2=2rs$, $\gamma (H)=1$,
and the canonical $nef$ element
$\wth=(-a,0,b)\mod \bz v$, if $\gamma (h)=1$ and the
conditions (i) and (ii) below hold:

(i) there exists $\wH\in N(Y)$ with
$\wH^2=2rs$, $\gamma(\wH)=1$, and such that
there exists an embedding
$$
f:K(\wH)\rightarrow K(h)
$$
of negative definite lattices such that
$$
\split
&f^\ast (K(h))=[K(\wH),2abcu^\ast (\wH)],\\
&f^\ast(u^\ast (h))+f^\ast(K(h))=m(a,b)c\,u^\ast(\wH)+f^\ast(K(h)).
\endsplit
\tag{5.1}
$$

(ii) For the $h$ and $f$ in (i), there exists a choice of $\pm$
such that
$$
\kappa(h)(k^\ast)=\pm \kappa(\wH)(f^\ast(k^\ast))
\tag{5.2}
$$
for any $k^\ast \in K(h)^\ast$.

The condition $\gamma (h)=1$ is necessary for any $Y$.

The conditions (i) and (ii) are necessary for a K3 surface $Y$ with
$\rho (Y)\le 19$ which is general for its Picard lattice $N(Y)$
in the following sense: the automorphism group of the transcendental periods
$(T(Y),H^{2,0}(Y))$ is $\pm 1$.

If $\rho (Y)=20$ the theorem is always valid, if $\gamma (h)=1$.
\endproclaim

Now assume that $\rk N(Y)=2$ and $h\in N(Y)$ has the degree $h^2=2ab$,
$\gamma (h)=1$. Applying Proposition 3.1.1 to $c=1$, we get

\proclaim{Proposition 5.2} Let $Y$ be a K3 surface with $\rk N(Y)=2$
equipped with a primitive $nef$ element $h$ of the degree
$h^2=2ab$, and $\gamma (h)=1$.

The pair $(N(Y),h)$ has the invariants $d\in \bn$ and
$\pm \nu\mod 2ab\in (\bz/2ab)^\ast$
such that $\nu^2\equiv d \mod 4ab$. We have: $\det N(Y)=-d$, and
$K(h)=h^\perp_{N(Y)}=\bz \delta_1$ where $\delta_1^2=-2abd$. Moreover,
$$
N(Y)=[h,\delta_1,\frac{(\nu h+\delta_1)}{2ab}],
\tag{5.3}
$$
$$
N(Y)=\{z=\frac{(xh+y\delta_1)}{2ab} \ |\ x,y\in \bz\ \text{and}\
x \equiv \nu y \mod 2ab \}.
\tag{5.4}
$$
We have $z^2=(x^2-dy^2)/2ab$.

For any primitive element $h^\prime \in N(Y)$ of the degree
$(h^\prime)^2 =h^2=2ab$ and the same invariant $\pm \nu$,
there exists an automorphism
$\phi\in O(N(Y))$ such that $\phi(h)=h^\prime$.
\endproclaim

From Proposition 3.1.1 applied to $(N(Y),\wH)$ where $\wH^2=2abc^2$ and
$\gamma(\wH)=1$, and Proposition 5.2 we get that
$\text{g.c.d}(d,2rs)=\text{g.c.d}(d,2ab)=1$, and we get
similar statement as Proposition 3.1.2.

\proclaim{Proposition 5.3}
Under conditions and notations of
Propositions 5.2, all elements
$\wH=(xh+y\delta_1)/(2ab) \in N(Y)$ with
$\wH^2=2abc^2$ are in one to one
correspondence with integral solutions $(x,y)$ of the equation
$x^2-dy^2=4a^2b^2c^2$ such that $x \equiv \nu y \mod 2ab$.

For $k\in \bn$ and $k>1$ the element $\wH$ is divisible by $k$
in $N(Y)$ and then it is not primitive, if and only if $k\,|\, (x,y)$ and
$x\equiv \nu y\mod 2abk$. Then $k^2|abc^2$. We call such a pair
$(x,y)$ as $\nu$-divisible by $k$. If $(x,y)$ is not $\nu$-divisible by
$k>1$, we call such a pair $(x,y)$ as $\nu$-primitive.
It is enough to consider only prime $k$. Thus, $\wH$ is primitive,
if and only if the corresponding pair $(x,y)$ is $\nu$-primitive.
If  $\text{g.c.d}(d,c)=1$ and $\wH$ is primitive, then $\gamma (\wH)=1$.
\endproclaim

The crucial statement is the statement which is similar to
Theorem 3.1.3:

\proclaim{Theorem 5.4}
Let $Y$ be a K3 surface with $\rho (Y)=2$, and $h\in N(Y)$ a $nef$ element
of the degree $h^2=2ab$ with $\gamma (h)=1$. As in Proposition 5.2,
we consider invariants $d$ and $\nu\in (\bz/2ab)^\ast$ of the
pair $(N(Y),h)$. Thus, $d\equiv \nu^2\mod 4ab$.
Assume that $\text{g.c.d}(c,d)=1$.

With notations of Proposition 5.2, all elements
$\wH=(xh+y\delta_1)/(2ab)\in N(Y)$ with $\wH^2=2abc^2$
satisfying Theorem 5.1 are in one to one correspondence with
integral solutions $(x,y)$ of the equation
$$
x^2-dy^2=4a^2b^2 c^2
\tag{5.5}
$$
which satisfy the conditions (i), (ii), (iii), (iv) below:

(i) $x\equiv \nu y\mod 2ab$;

(ii) $x\equiv \pm 2abc\mod d$;

(iii) $(x,y)$ belongs to one of $a$-series or $b$-series of solutions
defined below:

$a$-series: $b|(x,y)$ and $x-\nu y\equiv 0\mod 2a^2b$ where
$\nu^2\equiv d\mod 4a^2b$,

$b$-series: $a|(x,y)$ and $x-\nu y\equiv 0\mod 2ab^2$ where
$\nu^2\equiv d\mod 4ab^2$;

(iv) the pair $(x,y)$ is $\nu$-primitive: there does not exist
a prime $l$ such that $l|(x,y)$ and $x-\nu y\equiv 0\mod 2ab\,l$
(we always have $l^2|abc^2$ for such $l$), i. e. \newline
$\text{g.c.d}\left(x,y,(x-\nu y)/(2ab)\right)=1$.

In particular, by Theorem 5.1, for a general $Y$ with $\rho(Y)=2$,
the pair $(Y,h)$ is the moduli space $(\M,\wth)$ of coherent sheaves
over $Y$ with the Mukai vector $v=(r,H,s)$, where $H^2=2rs$,
$\gamma (H)=1$, and the canonical $nef$ element $\wth=(-a,0,b)\mod \bz v$,
if $\text{g.c.d}(d,c)=1$ and the equation $x^2-dy^2=4a^2b^2c^2$ has
an integral solution $(x,y)$ satisfying conditions (i)---(iv) above.
Moreover, a $nef$ primitive element $\wH=(xh+y\delta_1)/(2ab)\in N(Y)$
with $\wH^2=2abc^2$ defines the Mukai vector $(r,\wH,s)$ with the
moduli space $(\M,\wth)$ of coherent sheaves on $Y$
which is isomorphic to $(Y,h)$, if and only if $(x,y)$
satisfies the conditions (ii) and (iii) (it satisfies the
conditions (i) and (ii) since it corresponds to a
primitive element of $N(Y)$).
\endproclaim

We can deduce from Theorem 5.4 a similar result as Theorem 3.1.4

\proclaim{Theorem 5.5}
With conditions and notations of Theorem 5.4,
for a general $Y$ with $\rho (Y)=2$ and $\gamma(h)=1$,
we have $(Y,h)$ is the moduli space $(\M,\wth)$
of sheaves over $Y$ with the Mukai vector 
$v=(r,\wH,s)$ where $\wH^2=2rs$, $\gamma (\wH)=1$, 
and the canonical $nef$
element $\wth=(-a,0,b)\mod \bz v$, if and only if $\text{g.c.d}(c,d)=1$, 
for one of $\alpha=\pm 1$ there exists an integral solution $(p,q)$ 
for one of the $a$-series or the $b$-series:

The $a$-series:
$$
\cases
p^2-dq^2=4ac/\alpha \\
p\equiv \nu q \mod 2a\\
\text{g.c.d}(l_1,p,q)=1\ \text{for any prime }\ l_1|ac\ \text{and}\
l_1\not|b\\
p\not\equiv \nu q\mod 2al_2\ \text{for any prime}\  l_2|b\\
\endcases .
\tag{5.6}
$$
A solution $(p,q)$ of \thetag{5.6} gives a solution $(x,y)$
of Theorem 5.4 as the associated solution
$$
(x,y)=\pm (2abc+\alpha bdq^2,\ \alpha  bpq).
\tag{5.7}
$$

The $b$-series:
$$
\cases
p^2-dq^2=4bc/\alpha \\
p\equiv \nu q \mod 2b\\
\text{g.c.d}(l_1,p,q)=1\ \text{for any prime }\ l_1|bc\ \text{and}\
l_1\not|a\\
p\not\equiv \nu q\mod 2bl_2\ \text{for any prime}\  l_2|a\\
\endcases .
\tag{5.8}
$$
A solution $(p,q)$ of \thetag{5.8} gives a solution $(x,y)$
of Theorem 5.4 as the associated solution
$$
(x,y)=\pm (2abc+\alpha adq^2,\ \alpha apq).
\tag{5.9}
$$
All $\wH$ (or solutions $(x,y)$) of Theorem 5.4 are given by
\thetag{5.7} and \thetag{5.9} as associated solutions.
\endproclaim

Like in Theorems 3.1.5 and 3.1.6 we can interpret $(p,q)$ of Theorem 5.5
as elements $h_1=b(ph+q\delta_1)/(2ab)\in N(Y)$ for the $a$-series
and $h_1=a(ph+q\delta_1)/(2ab)\in N(Y)$ for the $b$-series.  We formulate
only the statement which is analogous to Theorem 3.1.6. It contains
the statement which is similar to Theorem 3.1.5.

\proclaim{Theorem 5.6}
Let $Y$ be a K3 surface with a $nef$ element $h$ such that
$h^2=2ab$ and $\gamma (h)=1$. Then $(Y,h)$ is
the moduli space $(\M,\wth)$ of sheaves over $Y$ with the Mukai vector
$v=(r,\wH,s)$ where $\wH^2=2rs$, $\gamma (\wH)=1$, and the canonical
$nef$ element $\wth=(-a,0,b)\mod \bz v$, if at least for one of
signs $\pm$ (i. e. $\alpha=\pm 1$) there exists $h_1\in N(Y)$ such
that elements $h$, $h_1$ are contained in a 2-dimensional sublattice
$N\subset N(Y)$ with $h\cdot N=\bz$, $\text{g.c.d}(\det (N),c)=1$,
and $h_1$ belongs to one of the $a$-series or the $b$-series:

The $a$-series:
$$
h_1^2=\pm 2bc,\  h\cdot h_1\equiv 0\mod b,\
h\cdot h_1\not\equiv 0\mod bl_1,\
h_1/l_2\not\in  N(Y)
\tag{5.10}
$$
where $l_1$ and $l_2$ are any primes such that $l_1|ac$ and
$l_1\not|b$, and $l_2|b$.

The $b$-series:
$$
h_1^2=\pm 2ac,\  h\cdot h_1\equiv 0\mod a,\
h\cdot h_1\not\equiv 0\mod al_1,\
h_1/l_2\not\in  N(Y)
\tag{5.11}
$$
where $l_1$ and $l_2$ are any primes such that $l_1|bc$ and
$l_1\not|a$, and $l_2|a$.

The conditions above are necessary, if
$\rho(Y)\le 2$ and $Y$ is a general K3 surface with
its Picard lattice, i. e. the automorphism group of the transcendental
periods $(T(Y), H^{2,0}(Y))$ of $Y$ is $\pm 1$.

\endproclaim

\remark{Important Remark} By Theorem 5.4 and the formulae \thetag{5.7}
and \thetag{5.9} we get
$$
\pm w(\wH)=
\cases
&-ch +\frac{\alpha (h\cdot h_1)h_1}{b},\ \text{if $h_1$
is from $a$-series,}\\
&-ch + \frac{\alpha (h\cdot h_1)h_1}{a},\ \text{if $h_1$
is from $b$-series}
\endcases
\tag{5.12}
$$
where $w\in W^{(-2)}(N(Y))$.
\endremark

\smallpagebreak

Using results of Sects 3 and 5 both, one can construct many
compositions of the correspondences from Sects 3 and 5
of the K3 surfaces $X\cong Y$ with itself.

\enddocument

\end